\documentclass[10pt,a4paper]{elsarticle}
\journal{Discrete Applied Mathematics}
\usepackage[latin1]{inputenc}
\usepackage{amsmath}
\usepackage{amsfonts}
\usepackage{amssymb}
\usepackage{wasysym}
\usepackage{graphicx}
\usepackage{geometry}
\usepackage{amsmath}
\usepackage{amssymb}
\usepackage{amsthm}
\usepackage{color}
\usepackage{tikz}
\newtheorem{mydef}{Definition}
\newtheorem{mylem}{Lemma}
\newtheorem{mythm}{Theorem}

\begin{document}

\begin{frontmatter}
\title{Implicit representation conjecture for semi-algebraic graphs}

\author{Matthew Fitch}
\address{Mathematics Institute, Zeeman Building,
University of Warwick,
Coventry CV4 7AL, United Kingdom \\ Email address: M.H.D.Fitch@warwick.ac.uk}

\begin{abstract}
The implicit representation conjecture concerns hereditary families of graphs. Given a graph in such a family, we want to assign some string of bits to each vertex in such a way that we can recover the information about whether 2 vertices are connected or not using only the 2 strings of bits associated with those two vertices. We then want to minimise the length of this string. The conjecture states that if the family is hereditary and small enough (it only has $2^{O(n\ln(n))}$ graphs of size $n$), then $O(\ln(n))$ bits per vertex should be sufficient. The trivial bounds on this problem are that: (1) some families require at least $\log_2(n)$ bits per vertex ; (2) $(n-1)/2+\log_2(n)$ bits per vertex are sufficient for all families.\\
\\
In this paper, we will be talking about a special case of the implicit representation conjecture, where the family is semi-algebraic (which roughly means that the vertices are points in some euclidean space, and the edges are defined geometrically, or according to some polynomials). We will first prove that the `obvious' way of storing the information, where we store an approximation of the coordinates of each vertex, does not work.  Then we will come up with a way of storing the information that requires $O(n^{1-\epsilon})$ bits per vertex, where $\epsilon$ is some small constant depending only on the family. This is a slight improvement over the trivial bound, but is still a long way from proving the conjecture.

\end{abstract}

\end{frontmatter}

\section{Introduction}

\begin{mydef}
Given a Euclidean space $\mathcal{S}$ and a finite set of symmetric polynomial eqalities and inequalities (which can be strict or not) on $\mathcal{S} \times \mathcal{S}$, the associated semi-algebraic family of graphs is the set of all graphs whose vertices are points in $\mathcal{S}$ and whose edges are exactly those pairs of points that satisfy all the constraints. 

\end{mydef}

For example, the family of closed disk graphs consists of all graphs whose vertices are closed disks in the plane and where edges indicate that two disks intersect. The vertices can be viewed as points in $\mathbb{R}^3$: $(x,y,r)$ where $(x,y)$ are the coordinates of the center of the disk and $r$ is the radius. There is an edge between $(x_1,y_1,r_1)$ and $(x_2,y_2,r_2)$ if and only if $(x_1-x_2)^2+(y_1-y_2)^2\leq (r_1+r_2)^2$.\\
\\
The family of unit disk graphs is defined in a similar way except all the radii are 1.\\
\\
Essentially, semi-algebraic families of graphs are made up of graphs that are defined geometrically or algebraically. They are very useful in graph theory because they are a good way of constructing graphs with certain properties. Work on semi-algebraic graphs has mostly been focused on specific families, such as the aforementioned family of disk graphs, which has applications in computational geometry \cite{breu1998unit}. However, there are a few general results. In 2005, Alon, Pach, Pinchasi, Radoi\v{c}i$\acute{c}^c$ and Sharir proved that given any semi-algebraic family of graphs, that every graph in it with $n$ vertices contains two subsets of vertices of size $\epsilon n$ (where $\epsilon$ is a constant), such that either all edges between them exist or none of them do. They also proved that there exists either a complete subgraph of size $n^\delta$ or an induced empty subgraph of size $n^\delta$ (where $\delta$ is another constant). \cite{alon2005crossing}\\
\\
In 2013, Blagojevi$\acute{c}$, Bukh and Karasev looked at algebraic methods while trying to solve the Tur\'{a}n problem for the complete bipartite graph $K_{s,s}$ and showed that one particular `natural' type of semi-algebraic graph cannot be used to construct a $K_{s,s}$-free graph with $\Theta(n^{2-1/s})$ edges \cite{blagojevic2013turan}.
\\
\\
The problem we are trying to solve in this paper is a special case of the Implicit Representation conjecture, first posed by Kannan, Naor and Rudich in 1992 \cite{kannan1992implicat}, which was also asked by Spinrad in 2003 \cite{spinrad2003efficient}. We want to come up with a method for storing graphs using the least number of bits per vertex. A hereditary family of graphs is one in which all induced subgraphs of every graph in the family are also in the family. Given such a hereditary family of graphs $\mathcal{G}$ (for example: disk graphs or unit disk graphs), let $\mathcal{G}^{(n)}$ mean the subset of graphs who have exactly $n$ vertices. For every $n$, we want a function $F^{(n)}: \mathcal{G}^{(n)} \rightarrow  [2^m]^n$, and a symmetric function $G^{(n)}: [2^m] \times [2^m] \rightarrow \{0,1\}$ such that for every graph $H\in \mathcal{G}^{(n)}$ and every pair of vertices $i,j$ in $H$, we have $G(F(H)_i,F(H)_j)=1$ if and only if there is an edge between $i$ and $j$. Furthermore, we want to minimise $m=m(n)$, which is the amount of information per vertex. The Implicit Representation Conjectures states that if there exists a constant $c$ such that the family $\mathcal{G}$ contains less than $2^{cn\ln(n)}$ graphs of size $n$ for all $n$, then there exists a constant $c'$ such that $m=c'\log_2(n)$ will be sufficient for every graph of size $n$ in the family. The Implicit Representation Conjecture has been proved for a large number of families by Atminas, Collins, Lozin and Zamaraev in \cite{atminas2015implicit}.\\
\\
A corollary of Warren's theorem (1968) \cite{warren} \cite[p.1763]{handbook} shows that semi-algebraic families do indeed have at most $2^{O(n\ln(n))}$ graphs of size $n$, so they do satisfy the conditions for the Implicit Representation Conjecture. We'll see the derivation of this corollary at the end of section 2.\\
\\
The trivial lower bound for this problem matches the conjecture, at $m=\log_2(n)$ since that is the amount of information required to identify a vertex amongst $n$. More specifically, if we use less than $\log_2(n)$ bits, then there are less than $n$ possible options for what the data can be, so by the pidgeonhole principle, there exist two vertices $i$ and $j$ with $F(H)_i=F(H)_j$. This means that their neighbourhoods are identical. However, if we let the graph be a path, then every vertex has a different neighbourhood, which is a contradiction. If the family is defined by the intersection of bounded non-trivial shapes (such as the disk graph), then it is fairly easy to see that we can draw a path using these shapes.\\
\\
A trivial upper bound is $m=\left\lceil\frac{n-1}{2}\right\rceil +\left\lceil\log_2(n)\right\rceil$. To achieve this, we write the vertices as $0,1,...,n-1$ in $\mathbb{Z}/n\mathbb{Z}$, and store this information using $\left\lceil\log_2(n)\right\rceil$ bits. Then for every vertex $i$, let $F(H)_i$ be a list of $\left\lceil\frac{n-1}{2}\right\rceil$ 0s and 1s, with a 1 in the $k$th position if and only if there is an edge between $i$ and $i+k$. For every pair $i$ and $j$, $G$ will then output the $(j-i)$th entry of $F(H)_i$ if $j-i$ is between 1 and $\left\lceil\frac{n-1}{2}\right\rceil$ and otherwise it will output the $(i-j)$th entry of $F(H)_j$. \\
\\
A natural idea we could have for semi-algebraic graphs would be to store integer approximations of the coordinates of all the vertices. This looks like a good idea because it is easy to store integers, and because the function $G^{(n)}$ is easy to compute (just evaluate all the polynomial inequalities). Unfortunately this does not work. In 2013, McDiarmid and M\"{u}ller \cite{mcdiarmid2013integer} proved that there exist unit disk graphs with $n$ vertices, every planar representation of which had to have four vertices $a,b,c,d$ for which $\frac{|a-b|}{|c-d|}>2^{2^{\Omega(n)}}$. If $a,b,c,d$ had integer coordinates, then one of $a$ or $b$ has to have a coordinate of size at least $2^{2^{\Omega(n)}}$. This requires $2^{\Omega(n)}$ bits to store which is even more than the trivial bound.\\
\\
In 2012, Kang and M\"{u}ller \cite{kang2012sphere} improved upon this result in two ways, firstly by showing that the dimension $k$ of the ambient space can be arbitrary, and secondly by replacing the integer approximations by rational approximations. They showed that for any $k\geq 2$, there exist unit $k$-ball graphs with $n$ vertices but for which every realisation of it in $\mathbb{R}^k$ had to have four vertices $a,b,c,d$ for which $\frac{|a-b|}{|c-d|}>2^{2^{\Omega(n)}}$. (A $k$-ball graph is defined the same as a disk graph except that the ambient space is of dimension $k$ instead of 2.) Then if $a,b,c,d$ had rational coordinates, then one of these four points has to have a coordinate with numerators or denominators of size at least $\sqrt[4]{2^{2^{\Omega(n)}}} = 2^{2^{\Omega(n)-2}}$. This requires $2^{\Omega(n)-2}$ bits to store, which is even more than the trivial upper bound. Thus, storing rational approximations of the coordinates of all the vertices does not work in general.\\
\\
In the first part, we will go even further, and ask whether we can store the coordinates as algebraic numbers instead of rational numbers or integers. However, this runs into the same problems, as we'll see shortly. \\
\\
In our second part, we find a very minor improvement on the upper bound that does work. It uses a result by Yao and Yao \cite{yao1985general}, and ideas about semi-algebraic sets from Alon et al. \cite{alon2005crossing}.

\begin{mythm}
Given a semi-algebraic family of graphs $\mathcal{G}$, there exists some constant $\epsilon>0$ such that we can store every graph of size $n$ in the family using $O(n^{1-\epsilon})$ bits per vertex for $n$ sufficiently large.\\
\\
 More specifically, for $m=O(n^{1-\epsilon})$ there exists a series of functions $F^{(n)}: \mathcal{G}^{(n)} \rightarrow  [2^m]^n$, and a symmetric function $G^{(n)}: [2^m] \times [2^m] \rightarrow \{0,1\}$ such that for every graph $H\in \mathcal{G}^{(n)}$ and every pair of vertices $i,j$ in $H$, we have $G(F(H)_i,F(H)_j)=1$ if and only if there is an edge between $i$ and $j$. 

\end{mythm}

\section{Semi-algebraic graphs}

\subsection{Simplification of the problem}

First of all, note that any equalities $f=g$ in the definition of a semi-algebraic family can be written as the combination of two inequalities $f\geq g$ and $f\leq g$, so that part of the definition was redundant. So we can assume without loss of generality that there are only inequalities.\\
\\
Secondly, we will note that we can reduce to the case where the semi-algebraic family is defined by only one equation. Indeed, if we have a semi-algebraic family defined by $k>1$ polynomials $f_1,f_2,...,f_k$ and we have a graph $G$ in this family with vertices $x_1,...,x_n\in \mathbb{R}^q$. Then the edge set of $G$ can be viewed as a boolean function of the edge sets of $k$ semi-algebraic graphs $G_1,G_2,...,G_k$, each with the same vertex set and defined by polynomials $f_1$, $f_2$, ... , $f_k$ respectively. If we can store each of these $G_i$s using $m$ bits per vertex, then by concatenation, we can store $G$ using $k \cdot m$ bits per vertex. So without loss of generality, we can assume that there is only a single polynomial inequality $f(x,y)\geq 0$ or $f(x,y)>0$ that defines the semi-algebraic family. \\
\\
Also note that the complement of a graph can be stored using the same number of bits as the original graph, by simply exchanging 0 and 1 in the output of the function $G$. Therefore we can without loss of generality assume the single polynomial inequality is of the form $f(x,y)\geq 0$.\\
\\
\\
Now suppose we have a semi-algebraic family of graphs whose vertices can be written as living in the space $\mathbb{R}^q$ and where for any $x,y\in \mathbb{R}^q$, $(x,y)$ is an edge if and only if $f(x,y)\geq 0$ (where $f$ is a polynomial). Let $d$ be the degree of $f$. Now for every vertex $x=(x_1,x_2,...,x_q)$ in  $\mathbb{R}^q$, we can replace it with a point $\tilde{x}$ consisting of all the terms of degree less than or equal to $d$, i.e., $\tilde{x}=(1,x_1,x_2,...,x_q,x_1^2,x_1x_2,x_1x_3,...,x_q^2,x_1^3,x_1^2x_2,...,x_q^3,...,x_q^d)$. This point exists in the space $\mathbb{R}^{\binom{q+d}{d}}$. Set $Q={\binom{q+d}{d}}$. We can then also rewrite the polynomial $f(x,y)$ as a bilinear function $f(x,y) = \tilde{x}^T M \tilde{y}$ where $M$ is a $Q \times Q$ matrix ($M$ is symmetric because $f$ was symmetric). \\
\\
\\
Note that given a fixed $\tilde{y}$, the set of solutions to the equation $z^T M \tilde{y} \geq 0$  forms a half-space, so for every vertex, the set of vertices adjacent to it is exactly those in the half-space. \\
\\
At this point, we already have everything we need about semi-algebraic graphs to complete the proof; however, the $\epsilon$ that we will get in the final result will be a function of the dimension $Q$ so decreasing $Q$ will improve the result slightly. So there is one more thing we can do: it is a standard property of bilinear forms that we can diagonalise them, and furthermore, we can make it such that there are only $1$s, $-1$s and $0$s on the diagonal. So without loss of generality, we can assume that

\[M=\left(\begin{matrix}
 1 & 0 & 0 & ... & 0 & 0 & ... & 0 \\
 0 & 1 & 0 & ... & 0 & 0 & ... & 0 \\
 0 & 0 & 1 & ... & 0 & 0 & ... & 0 \\\
 ... & ... & ... & ... & ... & ... & ... \\
 0 & 0 & 0 & ... & -1 & 0 & ... & 0 \\
  0 & 0 & 0 & ... & 0 & -1 & ... & 0 \\
 ... & ... & ... & ... & ... & ... & ... & ... \\
 0 & 0 & 0 & ... & 0 & 0 & ... & 0

\end{matrix}\right)\]

We can delete those coordinates for which $M$ has a zero on the diagonal because they do not impact the result. So without loss of generality, $M$ is a diagonal matrix with only $1$s and $-1$s on the 
diagonal. In particular, there are only $Q+1$ types of matrix in dimension $Q$. So solving the problem for just these few special cases is enough. Every other semi-algebraic family is a combination of matrices of this type after a change of basis. \\
\\

\textit{Example 1:\\}

For the unit disk graph in the plane, every vertex can be identified with its center: $(x,y)$. Then two disks $(x_1,y_1)$ and $(x_2,y_2)$ intersect if and only if $(x_1-x_2)^2+(y_1-y_2)^2\leq 4$. There is only a single inequality, so if we put this in bilinear form, we get a single matrix:

\[(x_1^2,x_1y_1,y_1^2,x_1,y_1,1)\left(\begin{matrix} 
0 & 0 & 0 & 0 & 0 & -1\\
0 & 0 & 0 & 0 & 0 & 0  \\
0 & 0 & 0 & 0 & 0 & -1\\
0 & 0 & 0 & 2 & 0 & 0\\
0 & 0 & 0 & 0 & 2 & 0\\
-1 & 0 & -1 & 0 & 0 & 4
\end{matrix}\right)\left( \begin{matrix} 
x_2^2 \\
x_2y_2 \\
y_2^2 \\
x_2 \\
y_2 \\
1
\end{matrix}\right) \geq 0\]

We can use a change of basis and then delete irrelevant coordinates to turn this matrix into:

\[\mathbf{x}'^T\left(\begin{matrix} 
1 & 0 & 0 & 0 \\
0 & 1 & 0 & 0 \\
0 & 0 & 1 & 0 \\
0 & 0 & 0 & -1 \\
\end{matrix}\right)\mathbf{y}' \geq 0\]

Note that the dimension of this matrix is $Q=4$.\\

\textit{Example 2:\\}

For the disk graph in the plane, every vertex can be identified with its center and its radius: $(x,y,r)$. Then two disks $(x_1,y_1,r_1)$ and $(x_2,y_2,r_2)$ intersect if and only if $(x_1-x_2)^2+(y_1-y_2)^2\leq (r_1+r_2)^2$. If we put this in bilinear form, we get

\[(x_1^2,y_1^2,r_1^2,x_1,y_1,r_1,1)\left(\begin{matrix} 
0 & 0 & 0 & 0 & 0 & 0 & -1\\
0 & 0 & 0 & 0 & 0 & 0 & -1\\
0 & 0 & 0 & 0 & 0 & 0 & 1\\
0 & 0 & 0 & 2 & 0 & 0 & 0\\
0 & 0 & 0 & 0 & 2 & 0 & 0\\
0 & 0 & 0 & 0 & 0 & 2 & 0\\
-1 & -1 & 1 & 0 & 0 & 0 & 0\\
\end{matrix}\right)\left( \begin{matrix} 
x_2^2 \\
y_2^2 \\
r_2^2 \\
x_2 \\
y_2 \\
r_2 \\
1
\end{matrix}\right) \geq 0\]

Using the change of basis, we can replace this by:

\[\mathbf{x}'^T\left(\begin{matrix} 
1 & 0 & 0 & 0 & 0 \\
0 & 1 & 0 & 0 & 0 \\
0 & 0 & 1 & 0 & 0 \\
0 & 0 & 0 & 1 & 0 \\
0 & 0 & 0 & 0 & -1 \\
\end{matrix}\right)\mathbf{y}' \geq 0\]

So in this case we have $Q=5$.\\
\\
\subsection{Proof that semi-algebraic families satisfy the conditions for the Implicit Representation Conjecture}

For this, we use Warren's theorem \cite{warren} \cite[p.1763]{handbook}: 

\begin{mythm}[Warren, 1968]
Suppose we have a set of $k$ real polynomials in $l$ variables of degree at most $d$ and $k\geq l$. If we split $\mathbb{R}^l$ into regions depending on the signs of all the polynomials (i.e., whether each polynomial is negative, positive or 0 at a given point in $\mathbb{R}^l$), we end up with at most $(8\mathbf{e}dk/l)^l$ regions. 
\end{mythm}

Suppose we have a semi-algebraic family $\mathcal{G}$, with associated euclidean space $\mathcal{S}$ and associated set of polynomial inequalities $\mathcal{P}$ and also let $n$ be an integer. We want to count the number of graphs in our family with $n$ vertices. A graph $G$ is in $\mathcal{G}$ if and only if there exist $n$ distinct points $x_1,x_2,...,x_n$ in $\mathcal{S}$, such that $\mathcal{P}(x_i,x_j)$ is true if and only if $(x_i,x_j)$ is an edge of $G$. So let $\mathcal{Q}_G$ be the set of polynomial inequalities: $\bigcup_{(i,j) \text{ edge} }\mathcal{P}(x_i,x_j) \; \; \cup \; \; \bigcup_{(i,j) \text{ non-edge} }\neg\mathcal{P}(x_i,x_j)$. Then $G$ is in  $\mathcal{G}$ if and only if there exist $x_1,x_2,...,x_n$ in $\mathcal{S}$ that satisfy $\mathcal{Q}_G(x_1,x_2,...,x_n)$. Notably, we can split $\mathcal{S}^n$ into regions depending on the signs of the polynomials of $\mathcal{Q}_G$, and then each region will have a unique graph associated with it.\\ 
\\
So how many regions are there? $\mathcal{Q}_G$ is a set of $\binom{n}{2}|\mathcal{P}|$ polynomial inequalities. The number of variables of these polynomials is $n\cdot\text{dim}(\mathcal{S})$, and the maximum degree is $d$. So by Warren's theorem, the number of regions is at most $\left(\frac{8\mathbf{e}d\binom{n}{2}|\mathcal{P}|}{n\cdot\text{dim}(\mathcal{S})}\right)^{n\cdot\text{dim}(\mathcal{S})}$ as long as $n$ is large enough. This is less than  $\left(4\mathbf{e}d|\mathcal{P}|n/\text{dim}(\mathcal{S})\right)^{n\cdot\text{dim}(\mathcal{S})} \leq 2^{cn\ln(n)}$ for some large enough constant $c$. This completes the proof and shows that semi-algebraic families do in fact satisfy the hypothesis of the Implicit Representation Conjecture.\\
\\
\\

\section{The `algebraic points' method does not work for disk graphs}

A natural way we can try to store disk graphs is to assume that the centers and radii of all the circles be algebraic and just store these numbers. However, this turns out to be worse than the trivial bound. This builds upon McDiarmid and M\"{u}ller \cite{mcdiarmid2013integer}, who prove that storing the centers and radii of all the circles as rational numbers does not work. As part of their construction, they use the paper by Goodman, Pollack and Sturmfels \cite{goodman1989coordinate}, where a planar configurations of lines requires exponentially many bits to store using integer coordinates.\\
\\
An important part of the proof is that there exists an infinite family of disk graphs such that for any disk representation of them, there are four centers $x$,$y$,$z$ and $t$ such that $\frac{|x-y|}{|z-t|}>2^{2^{\Omega(n)}}$. This family was constructed in \cite{mcdiarmid2013integer}. We claim that this family also requires at least $\Omega(n)$ bits if the centers are algebraic (and where we store algebraic numbers in the standard way: by storing their minimal polynomial as a list of integers).\\
\\
First, we will describe in more detail the standard way of storing algebraic numbers. If $x$ is an algebraic number, it has a minimal integer polynomial it is a solution to: $\sum_{i=0}^{i_{max}} a_i x^i =0$. Using the standard way for storing integers, $a_i$ uses $\Theta(\log_2(|a_i|+1))$ bits. Therefore storing the polynomial takes $\Theta(\sum_{i=0}^k \log_2(|a_i|+1))+\Theta(i_{max})$ bits. The polynomial also has $i_{max}$ solutions so we additionally need $\Theta( \log_2(i_{max}))$ bits to indicate which solution it is. Therefore it overall takes $\Theta(\sum_{i=0}^{i_{max}} \log_2(|a_i|+1))+\Theta(i_{max})$ bits to store an algebraic number in the standard way. We'll call this number $m(x)$.\\
\\
For simplicity, we'll consider the real plane on which our disk graph is drawn to be $\mathbb{C}$, so each center only requires a single algebraic number to describe it.\\
\\
Pick some integer $m$. What is the largest we can make $|x-y|$ given that $m(x)\leq m$ and $m(y)\leq m$? First of all, we know that $|x-y|\leq |x|+|y|$. Next, suppose that $\sum_{k=0}^{k_{max}}c_k x^k =0$ is the minimal polynomial of $x$ (so in particular, $c_0\neq 0$). Then if $|x|>\sum_{k=0}^{k_{max}-1} |c_{k}|$, we have $|c_{k_{max}} x^{k_{max}}| > \left(\sum_{k=0}^{k_{max}-1} |c_{k}|\right)|x|^{k_{max}-1}$. We also have $|x|>|c_0|\geq 1$ so we can continue with: $\left(\sum_{k=0}^{k_{max}-1} |c_{k}|\right)|x|^{k_{max}-1} \geq \sum_{k=0}^{k_{max}-1} |c_k x^k| \geq |\sum_{k=0}^{k_{max}-1} c_k x^k|$. Putting this string of inequalities together we get $|c_{k_{max}} x^{k_{max}}| > |\sum_{k=0}^{k_{max}-1} c_k x^k|$ which contradicts the fact that $\sum_{k=0}^{k_{max}}c_k x^k =0$. Therefore we must have $x\leq \sum_{k=0}^{k_{max}-1} |c_{k}|$. This is bounded above by $2^{ \log_2(\left(\sum_{k=0}^{k_{max}-1} |c_{k}|\right)} \leq 2^{\sum_{k=0}^{k_{max}-1} \log_2(c_k)} \leq 2^{m-1}$. Similarly, $y\leq 2^{m-1}$. So the we have an upper bound of $2^m$ for $|x-y|$. \\
\\
Now what is the smallest we can make $|x-y|$ for distinct $x$ and $y$? Say $x$ has minimal polynomial $\sum_{i=0}^{i_{max}} a_i x^i$ while $y$ has minimal polynomial $\sum_{j=0}^{j_{max}} b_j y^j$. Our aim will be to find a polynomial which has $a_{i_{max}}b_{j_{max}}(x-y)$ as its solution. To do so, we will express $(a_{i_{max}}b_{j_{max}}(x-y))^l$ as a linear combination of $\{x^i y^j \mid i<i_{max} \, ; \, j<j_{max}\}$ for all $l$ from 0 to $i_{max}j_{max}$. Then because the linear space spanned by $\{x^i y^j \mid i<i_{max} \, ; \, j<j_{max}\}$ has dimension $i_{max}j_{max}$ but we have $i_{max}j_{max}+1$ elements in the space, we know that there has to be a linear dependence between them, which is a polynomial with $a_{i_{max}}b_{j_{max}}(x-y)$ as its solution.\\

Pick an integer $l\leq i_{max}j_{max}$ and consider $(a_{i_{max}}b_{j_{max}}(x-y))^l$. We can develop it into $(a_{i_{max}}b_{j_{max}})^l\sum_{k=0}^l \binom{l}{k}x^k(-y)^{l-k}$, which is expression in terms of $\{x^i y^j|i\leq l\, ; \, j\leq l\}$. The sum of the absolute values of the coefficients is $|a_{i_{max}}b_{j_{max}}|^l \cdot 2^l$. For $k$ running from $k=l$ down to $k=i_{max}$, we can replace all instances of $x^k$ in this expression with $-\sum_{i=0}^{i_{max}-1}a_ix^{k-i_{max}+i}/a_ {i_{max}}$ by using the minimal polynomial for $x$. \\
Note that every time we do this step, the maximum exponent of $x$ drops by at least 1, so after $l-i_{max}+1$ steps, the maximum exponent of $x$ will be at most $i_{max}-1$. Also note that at the end, all coefficients will still be integers: although we divide by $a_ {i_{max}}$ each step, $(a_{i_{max}})^l$ divided all the starting coefficients and we only do $l-i_{max}+1\leq l$ steps. \\
What does this do to the sum of the absolute values of the coefficients? Well every time we do this operation, we multiply it by at most $\frac{\sum_{i=0}^{i_{max}-1}|a_i|}{|a_{i_{max}}|}$. We know that $\sum_{i=0}^{i_{max}-1}\log_2(|a_i|+1) \leq O(m)$ so $\sum_{i=0}^{i_{max}-1}|a_i| \leq 2^{O(m)}$ by concavity of the $\log_2$ function. Since we started with the sum of the absolute values of the coefficients at most $|a_{i_{max}}b_{j_{max}}|^l \cdot 2^l$ and we do this operation $l-i_{max}+1$ times, we end up with the sum of the absolute values of the coefs is at most $|b_{j_{max}}|^l \cdot |a_{i_{max}}|^{i_{max}-1}\cdot (2^{O(m)})^{l-{i_{max}}+1})\cdot 2^l$. Note also that $|a_{i_{max}}|\leq 2^{O(m)}$ so we end up with the sum of the absolute values of the coefs is at most $ 2^l\cdot|b_{j_{max}}|^l\cdot 2^{O(ml)}$\\
\\
We do the same operation with $y$, to end up with a linear formula for $(a_{i_{max}}b_{j_{max}})^l(x-y)^{l}$ in terms of $\{x^i y^j \mid i<i_{max} \, ; \, j<j_{max}\}$, and where the sum of all the absolute values of all the coefficients is at most $2^{O(ml)}$.\\
\\
Now if we do this for all $l$ between $0$ and $i_{max}j_{max}$, then we have $i_{max}j_{max}+1$ formulas inside the linear space generated by $\{x^i y^j \mid i<i_{max} \, ; \, j<j_{max}\}$. But this space has dimension $i_{max}j_{max}$, so our formulas must be lineally dependent. Remebering that each of our formulas represented some power of $a_{i_{max}}b_{j_{max}}(x-y)$, this linear dependence is equivalent to an integer polynomial of degree $i_{max}j_{max}$ that is zero when evaluated at $a_{i_{max}}b_{j_{max}}(x-y)$. Without loss of generality suppose that this polynomial is minimal; say it has degree $d$. We'll write this polynomial as $\mu(a_{i_{max}}b_{j_{max}}(x-y))^d = \sum_{k=0}^{d-1} \lambda_k (a_{i_{max}}b_{j_{max}}(x-y))^k$ where $\mu$ and all the $\lambda$s are integers. How big are the coefficients of this polynomial? \\
\\
We can work out what they are. Since the polynomial was chosen to be minimal, we know that the formulas for $(a_{i_{max}}b_{j_{max}}(x-y))^k$ , $k<d$, are all lineally independent. We can list all these formulas in an $i_{max}j_{max} \times d$ matrix of integers which we'll call $M$, where the rows are linearly independent:\\ 

\[\left(\begin{matrix}
1 \\
a_{i_{max}}b_{j_{max}}(x-y) \\
(a_{i_{max}}b_{j_{max}}(x-y))^2 \\
... \\
(a_{i_{max}}b_{j_{max}}(x-y))^{d-1} \\
\end{matrix}\right) = 
\left(\begin{matrix} 
1 & 0 & 0 & 0 & ... & 0 \\
0 & a_{i_{max}}b_{j_{max}} & -a_{i_{max}}b_{j_{max}} & 0 &  ... & 0 \\
0 & 0 & 0 & a_{i_{max}}^2b_{j_{max}}^2 & ... & 0 \\
. & . & . & . &  ... & . \\
. & . & . & . &  ... & . 
\end{matrix} \right)\left(\begin{matrix}
1 \\
x \\
y \\
x^2 \\
xy \\
y^2 \\
...\\
x^{i_{max}-1}y^{j_{max}-1}
\end{matrix}\right) \]
Meanwhile, we also have a similar formula for $(a_{i_{max}}b_{j_{max}}(x-y))^d$, which takes the form of a vector of integers of size $i_{max}j_{max}$. We'll call this vector $\textbf{v}$:

\[(a_{i_{max}}b_{j_{max}}(x-y))^{d} \; \; \; \; = \; \; \; \; \textbf{v}\; .\; \left(\begin{matrix}
1 \\
x \\
y \\
x^2 \\
xy \\
y^2 \\
...\\
x^{i_{max}-1}y^{j_{max}-1}
\end{matrix}\right)\]

This is a linear combination of the rows of $M$: $\forall l \; , \; \mu \textbf{v}_l = \sum_{k=0}^{d-1} \lambda_k M_{k,l}$. If we write let the vector of $\lambda_k$s be $\mathbf{\lambda}$ (of length $d$), this formula can be rewritten in vector and matrix form as: 

\[\mu \mathbf{v} = \mathbf{\lambda}M\]

Now pick some lineally independant subset $C$ of the columns of the matrix of size $d$. This gives us a $d \times d$ non-singular matrix $M'$ Let $\pi_C$ be the matrix of the orthogonal projection from the space generated by $\{x^i y^j \mid i<i_{max} \, ; \, j<j_{max}\}$ onto the space generated by $C$. Thus,  $M'= M \pi_C $. Also let $\textbf{v'}$ be the image of $\mathbf{v}$ via this projection, ie $\mathbf{v'}=\mathbf{v}\pi_C$. The above linear combination continues to hold after projection:

\[\mu \mathbf{v'} = \mathbf{\lambda}M'\]

But now we can find out exactly what our $\lambda_k$s are by simply using the equation: $\mathbf{\lambda} = \mu M'^{-1}\textbf{v'}$ (remember that $M'$ is non-singular).\\
\\
$M'$ in an integer matrix so $M^{-1}$ is a rational matrix. Moreover, each element of $M^{-1}$ can be found as the determinant of a minor of $M'$ divided by $\text{det}(M')$, so in particular, each element is a multiple of $1/\text{det}(M')$. Without loss of generality, we can set $\mu=\text{det}(M')$ so as to make $(\mu M^{-1})$ an integer matrix; this makes all the $\lambda_k$s be integers too. We know that the sum of the absolute values of all the coefficients in each row of $M'$ is at most $2^{O(ml)}$, so we get that $|\text{det}(M')|$ is at most $d!(2^{O(ml)})^d=2^{O(md^2)}$. Moreover, the determinant of each minor of $M'$ is also at most $2^{O(md^2)}$ in absolute value so the elements of the matrix $(\mu M^{-1})$ are also at most $2^{O(md^2)}$ in absolute value. We know from before that the sum of the absolute values of the coefficients in $\mathbf{v'}$ is bounded above by $2^{O(md)}$ so for every $k$, $\lambda_k=d\cdot 2^{O(md^2)}\cdot 2^{O(md)} = 2^{O(md^2)}$.  \\
\\
Putting all of this together (remembering that $d\leq i_{max}j_{max}$), we get an integer polynomial that is 0 at $a_{i_{max}}b_{j_{max}}(x-y)$, that has of degree at most $i_{max}j_{max}$, and where all the coefficients are at most $2^{O(m[i_{max}j_{max}]^2)}$. Now both $i_{max}$ and $j_{max}$ are $\leq O(m)$ so this means the polynomial has degree at most $O(m^2)$ with coefficients at most $2^{O(m^5)}$. Say this polynomial is $\sum_{k=0}^{k_{max}} c_k (x-y)^k$.\\
\\
If we assume that $|x-y|<\frac{1}{\sum_{k=1}^{k_{max}} |c_k|}$, then $|\sum_{k=1}^{k_{max}} c_k (x-y)^k| \leq \left[\sum_{k=1}^{k_{max}} |c_k|\right]|x-y| < 1 \leq |c_0|$, which contradicts the polynomial being 0.  Therefore $|x-y|\geq\frac{1}{\sum_{k=1}^{k_{max}} |c_k|} = \frac{1}{O(m^2)}2^{\Omega(-m^5)} = 2^{\Omega(-m^5)}$\\
\\
\\
The ratio between the smallest possible value of $|x-y|$ and the largest is thus of order $2^{O(m^5)} \cdot 2^{O(m)} = 2^{O(m^5)}$. When we use the  special graph whose largest ratio is always at least $2^{2^{\Omega(n)}}$, we get that $m$ must be of order at least $2^{\Omega(n)/5}$. This is worse than our trivial upper bound of $m=(\frac{1}{2}+o(1))n$.

\section{A small improvement on the upper bound}

We will now present a method that allows us to store the information about a semi-algebraic graph using $O(n^{1-\epsilon})$ bits per vertex, where $\epsilon>0$ is a constant depending only on the semi-algebraic family. As a reminder, thanks to the work in section 2, we can assume without loss of generality that our semi-algebraic family is defined by only a single inequality: $xy$ is an edge of the graph if and only $f(x,y)\geq 0$ where $f$ is a symmetric polynomial function.\\
\\
We will start by doing the simpler case, where the vertices are in general position, that is to say, that $f(x,y)\neq 0$ for every pair of vertices $x$ and $y$. Then later we will do the harder case where we allow $f(x,y)=0$.

\subsection{The case where the vertices are in general position}

 We'll also assume for the moment that for every pair of vertices $(x,y)$, we never have $f(x,y)= 0$; we shall deal with that case at the end.\\
\\
In this case, for every vertex, the set of vertices adjacent to it is just a half-plane, which has a hyperplane as boundary. This is useful because it means we can use the following theorem:

\begin{mythm}[Yao and Yao, 1985 \cite{yao1985general}] Given a continuous and everywhere positive probability density function on $\mathbb{R}^Q$, there exists a partition of $\mathbb{R}^Q$ into $2^Q$ regions, each with mass equal to $1/2^Q$ such that every hyperplane in $\mathbb{R}^Q$ must avoid at least one of these regions. \\
\\
Moreover, these regions are convex polyhedral cones and all the cones have a common apex, called the center.  

\end{mythm}

A corollary of this theorem is the discrete version of it: 

\begin{mylem} Given a finite set $V$ of $n$ points in $\mathbb{R}^Q$, there exists a partition of $\mathbb{R}^Q$ into $2^Q$ regions, each of which contains between $\lfloor n/2^Q \rfloor$ and $\lfloor n/2^Q \rfloor + 2^Q-1$ of the points, such that every hyperplane in $\mathbb{R}^Q$ must avoid the interior of at least one region. \\
\\
Moreover, these regions are convex polyhedral cones and all the cones have a common apex, called the center.  
\end{mylem}

\textit{Proof of the lemma:}

Pick some small $\epsilon>0$. We replace each point $x\in V$ by a continuous density function on the ball of radius $\epsilon$ centered at $x$ whose total weight is $(1-\epsilon)/n$. We also add a continuous everywhere positive density function of total weight $\epsilon$. Adding up all of these densities gives a continuous everywhere positive probability density function on $\mathbb{R}^Q$, which means we can apply Yao and Yao's theorem. This splits the space into $2^Q$ convex polyhedral cones with a common apex, such that each has total weight $1/2^Q$. Let $A$ be the convex hull of the collection of balls centered at the points of $V$ of radius $\epsilon$ (i.e., $A$ is a bounded convex region of weight at least $1-\epsilon$).\\
\\
Suppose we are given an $\epsilon$, together with a polyhedral decomposition as in the lemma. For every vertex $x$ of $V$, we say it borders a certain region if the ball of radius $\epsilon$ centered around $x$ intersects the region. The information about which vertices border which regions will be called the configuration of the polyhedral decomposition. There are $2^Q$ regions so there are at most $2^{2^Q}$ possibilities for each vertex, and there are $n$ vertices so there are at most $n^{2^{2^Q}}$ possible configurations in total. This is finite therefore as $\epsilon \rightarrow 0$, there exists a configuration $C$ that occurs infinitely often. So we can pick a decreasing sequence of $\epsilon$s together with a corresponding collection of polyhedral cones in configuration $C$.\\
\\
Note that we do not have to use the Axiom of Choice to label the regions. Indeed, each region is a polyhedral cone and therefore has some defining equation. We can simply order these equations in lexicographic ordering and use that to label the regions. \\
\\
We'll define the \emph{shape} graph of one of these polyhedral decompositions as follows. The set of vertices is the collection of the faces of the decomposition (of any dimension). Two vertices are connected by an edge if one the corresponding faces contains the other. Since there are $2^Q$ regions, there are at most $2^{2^Q}$ faces, and therefore there are at most $2^{2^{2^{Q+1}}}$ possible shape graphs. This is finite so there is one particular shape graph that occurs infinitely often. So we pick a subsequence of $\epsilon$ such that all the polyhedral decompositions have the same shape graph. This means it now makes sense to talk about a particular face of the sequence of decompositions. \\
\\
Recall that the polyhedral decompositions that Yao and Yao's theorem construct have centers, that is to say, the point that is a common apex to all the regions. Moreover, these centers will stay inside $A$. To see why, assume not, and that we have a decomposition whose center is outside $A$. Pick a tangent hyperplane $T$ to $A$ that separates it from the center. Then every region of the polyhedral decomposition has weight at least $2^{-Q}$, while the exterior of $A$ has weight smaller than $\epsilon$, so when $\epsilon< 2^{-Q}$, every region's intersection with $A$ has to have non-zero weight. This implies in particular that every region's intersection with $A$ is non-empty. Since every region also has an apex at the center and is connected, that means every region has to cross $T$. That means $T$ is a hyperplane that fails to avoid a region, contradicting Yao and Yao's theorem. Therefore the center has to be within the bounded region $A$, so there is a subsequence of $\epsilon$s such that the centers converge to some point $M$.\\
\\
Now a given facet of a decomposition in our sequence is a part of a hyperplane that passes through the center of the decomposition. Moreover these centers of the decompositions converge to $M$ so the hyperplanes eventually have to pass within some small distance $\delta>0$ of $M$. Since the space of hyperplanes passing within $\delta$ of $M$ is compact, there is a hyperplane $H$ passing through $C$ and a subsequence of decompositions such that our given facet converges to a part of $H$. Repeat for all the other facets of the decomposition. \\
\\
If we take $M$ together with all the hyperplanes passing through it that we constructed and put the facets where they're supposed to be on the hyperplanes, we end up with a polyhedral decomposition of the space, which we'll call $P$. Moreover, our sequence of decompositions tends towards $P$ as $\epsilon \rightarrow 0$. We know that each decomposition in the sequence has configuration $C$. So for every vertex $x$ and for all $\epsilon>0$, $x$ is within $\epsilon$ of all the regions it is supposed to border according to $C$. Therefore $x$ is in the closure of all the regions of $P$ it is supposed to border. \\
 \\
$P$ therefore has the property that, for every region, its closure contains at least $n/2^Q$ points. In fact, we can go further and say that for any set $T$ of $t$ regions, the union of their closures contains at least $tn/2^Q$ points. Then by Hall's marriage theorem, there exists a way of associating disjoint sets of $\lfloor n/2^Q \rfloor$ points to each region such that the points are inside the closure of that region. There are at most $2^Q-1$ points left over, which we put in whichever region can accept them. \\
 \\
We have therefore created a partition of the points of $V$ into $2^Q$ regions, such that the closure of each region contains between $\lfloor n/2^Q \rfloor$ and $\lfloor n/2^Q \rfloor + 2^Q-1$ of the points, and such that every hyperplane in $\mathbb{R}^Q$ must avoid at least one of the interiors of a region. Thus the lemma is proved.
\begin{flushright}$\square$\end{flushright}

We will now use the discrete version of Yao and Yao's Theorem several times to store the information about a given semi-algebraic graph. This will mean we end up with a large number of partitions, each with their corresponding $2^Q$ regions. A rough outline of the method is as follows: given a vertex $x$, we will need the information about some of the regions it is in, and given a hyperplane perpendicular to a vertex $y$, we will need the information about some of the regions the hyperplane avoids. With this information, we will find a region that contains $x$ but that is avoided by the hyperplane corresponding to $y$, and this will tell us whether $x$ is adjacent to $y$ or not. Note that we will not need to store any specific information about the partition or the shape of the regions constructed in Yao and Yao's Theorem, only what vertices the regions contain and what hyperplanes they avoid. \\
\\
More precisely, the function $F^{(n)}$ that we contruct will do two things: Given a vertex $x$, it will provide an address $A(x)$ that will identify it by providing information about which regions of the Yao and Yao partitions it is in. Secondly, it will provide a tree structure $B(x)$ which defines which addresses it has an edge to and which ones it does not. This takes the form of a tree with labels on all its nodes. \\
\\
\textit{The address:} Apply the lemma to split the space into $2^Q$ regions, each of which contains at most $\lfloor n/2^Q \rfloor +2^Q-1$ vertices. We'll number these regions $1,2,...,2^Q$ and for each vertex $x$, we will then store the information about which region it is in as the first line of the address, and we call this $A_1(x)$. This takes $Q$ bits per vertex.\\
\\
Then repeat this process with every region, splitting each further into $2^Q$ subregions, then splitting each subregions into $2^Q$ subsubregions, etc. Continue until there are only $4^Q$ vertices in any given subregion. This will end in at most $s=\lceil \log_2\left(\frac{n-2^Q+1}{4^Q-2^{Q}+1}\right)/Q \rceil \cong \log_2(n)/Q-2$ steps. We then split this final subregion into its constituent points. Since there are at most $4^Q$ vertices in this subregion, this final decomposition also only takes $2Q$ bits. \\
\\
Thus, each vertex $x$ has a unique address $A(x)$ which takes the form of a string of $s+2$ numbers: $(i_1,i_2,...,i_{s+2})$ where each $i_w$ is an integer between $1$ and $2^Q$. The total amount of information stored in each vertex for the address ends up being $Q(s+2) \cong \log_2(n)$. \\

\textit{The tree-structure:}

Given a vertex $y$, we will construct the labelled tree $B(y)$ by induction. At step 0, we start with just the root node and leave it without a label. Throughout the construction, all the nodes in the tree can be matched onto certain partial addresses. The root node gets matched onto the empty label.\\
\\
Suppose we are at a certain step of the algorithm and that there exists an unlabeled leaf node in the tree. Say it can be matched to the partial address $(i_1,i_2,...,i_l)$. The first thing we do is give it $2^Q$ child nodes. We will match each of these child nodes to the addresses $(i_1,i_2,...,i_l,t)$ for every value of $t$ between 1 and $2^Q$. Now  because the graph is semi-algebraic, we know that there exists some half-space such that for every other vertex $x$, $x$ is connected to $y$ if and only if $x$ is in that half-space. This half-space has a hyperplane as its boundary, which we'll call $\mathcal{H}$.  Remember that when writing the address we split the region $(i_1,i_2,...,i_l)$ into $2^Q$ subregions using theorem 1 so theorem 1 tells us that $\mathcal{H}$ must avoid at least one of the interior's of a subregion. We will write a list of all the subregions whose interior it avoids on node $(i_1,i_2,...,i_l)$. Say it avoids the interior of the $t$th subregion. Now this subregion's interior is either entirely contained within the half-space or it is entirely disjoint from it. In other words, either all the vertices in the interior of the subregion are adjacent to $y$ or none of them are. In fact, because we assumed that $f(x,y)\neq 0$, for all $x,y$, $\mathcal{H}$ will not pass through any of the vertices, so this also extends to vertices on the boundary of the subregion. So we know that either all vertices in the region are adjacent to $y$ or none of them are. If all the vertices are adjacent, we will write a "1" on the $t$th child node. Otherwise write a "0" on the $t$th child node. Leave all the other child nodes unlabeled for the time being.\\
\\ 
Continue in this fashion until the only empty nodes in the tree correspond to sets of size less than $4^Q$. This will eventually happen at step number $s$. For each of these empty nodes, write down the size of the corresponding set on the node, and then creates $4^Q$ child nodes, each with a '1' or a '0' to indicate whether it is or is not adjacent to $y$. Thus, we will end up with a tree of depth at most $s+1$. This tree is comprised of some nodes with $2^Q$ child nodes; call these "splitting nodes" (except the final splitting nodes which have $4^Q$ children instead). The rest of the nodes just have a single number, "0" or "1" on them. We call these "leaf nodes". \\
\\
\textit{The function G:}
The function $G$ is simple to construct. Given two vertices $x$ and $y$, look at $x$'s address. Say it is $(i_1,i_2,...,i_{l+2})$. Now look at $y$'s tree. Travel through this tree by starting at the root node, and at every step $l$, if we are at a splitting node, then go to the $i_l$th child node. Eventually we will reach a leaf node and at that point, we should be able to read "1" or "0". If there is a "1", that means there is an edge between $x$ and $y$. If there is a "0", that means there is not.\\

\textit{Information used:}
What is the maximum amount of information required to store this tree? The structure (whether a certain node has a child or not) is entirely determined by the numbers written on each node, so we only have to count up the total information stored in the numbers. We'll work backwards from the end. \\
Each leaf node has either a "0" or a "1" so we have 1 bit per leaf node. \\
The final splitting nodes at the end have at most $4^Q$ children, each requiring 1 bit, so that's $4^Q$ bits for the children. It also stores how many children it has which requires an additional $2Q$ bits. So $4^Q+2Q$ bits suffice to store a splitting node at depth $l$ with all its descendants.\\
\\
Let $\alpha=8^Q-2\cdot4^Q+2\cdot2^Q-3Q-1$. We will prove by induction that $\frac{\alpha(2^Q-1)^m - (Q+1)}{2^Q-2}$ bits suffices to to store a splitting node at depth $s-m$ together with all its descendants. When $m=0$, it's easy to check our choice of $\alpha$ makes this hold.\\
\\
Now suppose we have a splitting node $\mathbf{i}$ which is at depth $ s-m$ for some $m\geq 1$. How much information suffices for it and all its descendants? Suppose it has $a$ leaf nodes adjacent. Each of these uses $1$ bit for itself, and another $Q$ bits to be put on the list of leaf nodes at $\mathbf{i}$, for a total of $a(Q+1)$ bits. The other $2^Q-a$ nodes are all splitting nodes, so by the induction hypothesis, each can be described using only $\frac{\alpha(2^Q-1)^{m-1} - (Q+1)}{2^Q-2}$ bits. Totalling everything up, we get: $\frac{\alpha(2^Q-1)^{m-1} - (Q+1)}{2^Q-2}(2^Q-a) + a(Q+1)$. Since $a\geq 1$, we get that this is less than: $\frac{\alpha(2^Q-1)^m}{2^Q-2} - \frac{(Q+1)(2^Q-1)}{2^Q-2} + (Q+1) = \frac{\alpha(2^Q-1)^m}{2^Q-2} -\frac{Q+1}{2^Q-2}$.\\
\\
Therefore by induction, each splitting node at depth $s-m$ together with all its descendants can be described using only $\frac{\alpha(2^Q-1)^m}{2^Q-2} -\frac{Q+1}{2^Q-2}$ bits. Therefore the total number of bits that suffices to store the entire tree is $\frac{\alpha(2^Q-1)^{l}}{2^Q-2} -\frac{Q+1}{2^Q-2}$\\
\textit{Summing it all up:\\}
Summing the contribution from the address $A(y)$ and the contribution from the tree $B(y)$, we get that the total maximum number of bits that suffices to store $F^{(n)}(y)$ is:\\
\begin{eqnarray*}
& & Q(s+2) + \frac{\alpha(2^Q-1)^{s}}{2^Q-2} -\frac{Q+1}{2^Q-2} \\
& = & \frac{\alpha(2^Q-1)^{\left(\frac{\log_2(n)}{Q}-2\right)}}{(2^Q-2)}(1+o(1)) \\
& = & \frac{\alpha}{(2^Q-2)(2^Q-1)^2}2^{\log_2(2^Q-1)\left(\frac{\log_2(n)}{Q}\right)}(1+o(1)) \\
& = & \frac{\alpha}{(2^Q-2)(2^Q-1)^2}n^{\log_2(2^Q-1)/Q}(1+o(1)) \\
& = & n^{(1-\frac{1}{Q2^Q})(1+o(1))}
\end{eqnarray*}

When $n$ is large, this is ever so slightly better than the trivial upper bound of $\lceil (n-1)/2 \rceil + \lceil\log_2(n)\rceil$.\\
\\
\\
\subsection{The case where $f(x,y)=0$}

The reason the previous method might not work in this case is that if $f(x,y)=0$, then $x$ will be on the hyperplane $\mathcal{H}$ corresponding to $y$. It is possible for $x$ to be on the boundary of its region, and that $\mathcal{H}$ is tangent to that region. Then the interior of region containing $x$ would be completely on one side of $\mathcal{H}$, so if we used that method, then we might get incorrect information about the edge $x,y$.\\
\\
The way we fix this is we will consider the boundaries of regions to be  full regions themselves, each of which will get their own addresses. However, the key thing to note here is that every boundary region will have smaller dimension. More formally, start with the closures of the $2^Q$ regions of the original decomposition. If two regions intersect, then their intersection gets subtracted from both of the original regions, and is instead counted as a region of its own. Repeat this process until there are no more intersections. Since there were $2^Q$ parts originally, there are at most $2^{2^Q}-1$ regions in the new decomposition. Also for every $d$ between $1$ and $Q$, there are at most $\binom{2^Q}{1+Q-d}$ regions of dimension $d$.\\
\\
\textit{Label:}

When we store the address $A(x)$ of a point $x$, we might need to write down some of these new boundary regions in the address if $x$ happens to be in one of them. However, we claim that the address can still be written using only $2\lceil\log_2(n)\rceil + 2Q^2$ bits. Indeed, every time we have a region of dimension $Q$ with $n$ points in it, we subdivide it into $2^Q$ subregions of dimension $Q$ that contain at most $n/2^Q$ points and for every $d$ between $1$ and $Q-1$, at most $\binom{2^Q}{1+Q-d}$ subregions of dimension $d$, each of which contains at most $n$ points. There is a single region of dimension 0: the center of the decomposition, which obviously contains at most 1 point. \\
\\
When $n\leq 2^Q$, then we can easily decompose using only $\lceil\log_2(n)\rceil$ bits, which is well within the bound (by a factor of 2). For $n$ is larger, there are 4 cases:\\
\\
\textbf{Case 1: subregions of dimension $Q$:} For points in the subregions of dimension $Q$, we use the induction hypothesis to say that the last part of the addresses can be written in $2\lceil\log_2(n/2^Q)\rceil +2Q^2$ bits. As there are $2^Q$ such subregions, we can indicate which one they are in using an additional $Q$ bits. Finally, we use $\lceil\log_2(Q)\rceil$ bits at the start to indicate what $d$ is. Therefore their full addresses can be written using $2\lceil\log_2(n/2^Q)\rceil+2Q^2+Q+\lceil\log_2(Q)\rceil = 2\lceil\log_2(n)\rceil-2Q+2Q^2+Q+\lceil\log_2(Q)\rceil = 2\lceil\log_2(n)\rceil+2Q^2-(Q-\lceil\log_2(Q)\rceil)$ bits. Since $Q\geq \lceil\log_2(Q)\rceil$, this works.\\
\\
\textbf{Case 2: subregions of dimension $d$, where $1\leq d\leq Q-1$: }By the induction hypothesis, the last part of the addresses of points in the subregions of dimension $d$ can be written using $2\lceil\log_2(n)\rceil + 2d^2$ bits. As there are at most $\binom{2^Q}{1+Q-d}$ such subregions, we can write identify which one it is using $Q(1+Q-d)$ bits. Finally, we use $\lceil\log_2(Q)\rceil$ bits to indicate what $d$ is. Therefore we can write the addresses of these points using $2\lceil\log_2(n)\rceil + 2d^2+Q(1+Q-d)+\lceil\log_2(Q)\rceil = 2\lceil\log_2(n)\rceil + 2Q^2 + (2d^2-Qd-Q^2+Q+\lceil\log_2(Q)\rceil)$ bits. The worst case scenario for $d$ is either when $d=1$ or when $d=Q-1$. \\

When $d=1$, we have the number of bits is $2\lceil\log_2(n)\rceil + 2Q^2 + (2-Q^2+\lceil\log_2(Q)\rceil)$. But now for this case to even appear, we need, $Q\geq 2$ so $Q^2\geq 2+\lceil\log_2(Q)\rceil)$ so this works.\\

When $d=Q-1$, we have the number of bits is $2\lceil\log_2(n)\rceil + 2Q^2 + (-2Q+2+\lceil\log_2(Q)\rceil)$. But as before, $Q\geq 2$ so $2Q-2\geq  \lceil\log_2(Q)\rceil)$ so this works.\\
\\
\textbf{Case 3: subregions of dimension 0: } This only happens if $x$ is directly at the center of the decomposition. Then we do not need any extra information to identify $x$. We only need $\lceil \log_2(Q) \rceil $ bits to indicate that $d=0$. This is well below $2\lceil\log_2(n)\rceil+2Q^2$ so this easily works.\\
\\
Therefore, by induction we can write down the new address of every vertex $x$ using $2\lceil\log_2(n)\rceil+2Q^2$ bits regardless of what subregions $x$ is in.\\
\\
\\
\textit{Tree Structure:} We also need to remake the tree structure $B(y)$ using similar methods. Instead of $2^Q-1$ branches at every splitting node, we'll end up with $1+\sum_{d=1}^{Q} \binom{2^Q}{1+Q-d}$ branches. However, all the new branches will have far less information on them, with the end result that the entire tree does not require that much more information to store. More precisely, we can store the tree for $n$ vertices in dimension $Q$ in $F(Q,n)=c_Q(2^Q-1)^{\log_2(n)/Q} \, - \, c'_Q(2^{Q-1}-1)^{\log_2(n)/{(Q-1)}}$ bits for some sufficiently large $c'_Q$ and $c_Q$.  \\
\\
For the subregions of dimension $Q$, we do a similar thing to last time . We know that there is at least one subregion of dimension $Q$ that avoids the hyperplane $\mathcal{H}$ coreesponding to $y$, which means all the elements of this subregion are either all adjacent or all non-adjacent from $y$. We'll say that there are in fact $a\geq 1$ subregions of dimension $Q$ that avoid $\mathcal{H}$. We can identify each one using $Q$ bits, and then we need to add 1 more bit to say whether all elements are adjacent or non-adjacent to $y$. For all the others, we know that since each contains at most $n/2^Q$ points, by the induction hypothesis that can store all the information using $F(Q,n/2^Q)$ bits. In total, we can store this information using $a(Q+1) + (2^Q-a)F(Q,n/2^Q)$ bits. Since $F(Q,n/2^Q)>Q+1$, the worst case scenario is when $a=1$.\\
\\
For the subregions of smaller dimension, we know that there are less than $2^{2^Q}$ of them. Therefore we can identify each one using $2^Q$ bits. Each of these also needs its internal information storing. By the induction hypotesis, we can store each one using $F(d,n)$ bits where $d$ is its dimension. Since $d\leq Q-1$, we we know that this is less than $F(Q-1,n)$. In total, we can store all the information about these subregions using $2^{2^Q}(2^Q+F(Q-1,n))$ bits. Adding this all up, we get that we can store all the information about a region of dimension $Q$ with $n$ points in it using information:

\begin{eqnarray*}
& &(Q+1) + (2^Q-1)F(Q,n/2^Q) + 2^{2^Q}(2^Q+F(Q-1,n)) \\
& = & (Q+1) + c_Q(2^Q-1)(2^Q-1)^{\log_2(n/2^Q)/Q}- c'_Q(2^Q-1)(2^{Q-1}-1)^{\log_2(n/2^Q)/(Q-1)} + 2^{2^Q+Q}\\
& & + c_{Q-1}2^{2^Q}(2^{Q-1}-1)^{\log_2(n)/(Q-1)} - c'_{Q-1}2^{2^Q}(2^{Q-2}-1)^{\log_2(n)/(Q-2)}\\
&\leq &  (Q+1 + 2^{2^Q+Q}) + c_Q(2^Q-1)^{\log_2(n)/Q} + \left[c_{Q-1}2^{2^Q}-c'_Q\frac{(2^Q-1)}{(2^{Q-1}-1)^{Q/{(Q-1)}}}\right](2^{Q-1}-1)^{\log_2(n)/(Q-1)}\\
& \leq &  c_Q(2^Q-1)^{\log_2(n)/Q} + \left[\left(Q+1+2^{2^Q+Q}+c_{Q-1}2^{2^Q}\right)-c'_Q\frac{(2^Q-1)}{(2^{Q-1}-1)^{Q/{(Q-1)}}}\right](2^{Q-1}-1)^{\log_2(n)/(Q-1)}
\end{eqnarray*}

If this ends up being less than or equal to $c_Q(2^Q-1)^{\log_2(n)/Q} \, - \, c'_Q(2^{Q-1}-1)^{\log_2(n)/{Q-1}}$, then the induction would be complete. This is equivalent to:

\[\left[\left(Q+1+2^{2^Q+Q}+c_{Q-1}2^{2^Q}\right)-c'_Q\frac{(2^Q-1)}{(2^{Q-1}-1)^{Q/{(Q-1)}}}\right] \leq -c'_Q\]

\noindent which is again equivalent to 

\[c'_Q\left(\frac{(2^Q-1)}{(2^{Q-1}-1)^{Q/{(Q-1)}}}-1\right)\geq Q+1+2^{2^Q+Q}+c_{Q-1}2^{2^Q}.\]

But now $\frac{(2^Q-1)}{(2^{Q-1}-1)^{Q/{(Q-1)}}}> 1$ because $\frac{(2^Q-1)^{Q-1}}{(2^{Q-1}-1)^{Q}}> 1$ because $(2^Q-1)^{Q-1} > (2^{Q-1}-1)^{Q}$. Therefore we can pick a $c'_{Q}=\frac{Q+1+2^{2^Q+Q}+c_{Q-1}2^{2^Q}}{\frac{(2^Q-1)}{(2^{Q-1}-1)^{Q/{(Q-1)}}}-1}$ and it will work.\\
\\
Then we can pick $c_Q$ large enough to make the initial conditions of the induction hold (i.e., when $n< 2^Q$, we can describe everything using $F(Q,n)$ bits). If we do this too, then by induction, we'll get for all $n$, we can describe everything using $F(Q,n)$ bits.\\
\\
Then by induction on $Q$, we'll have a series of numbers $c_Q$ and $c'_Q$ such that we can describe the entire tree using $c_Q(2^Q-1)^{\log_2(n)/Q} \, - \, c'_Q(2^{Q-1}-1)^{\log_2(n)/{(Q-1)}}$ bits. \\
\\
This completes the proof. We have a method of storing the information about the edges which uses $O(n^{\log_2(2^Q-1)/Q})$ bits for every vertex.\\

\textit{Examples:\\}

 For the category of unit disk graphs, we have $Q=4$ which means this takes $O(n^{0.976723})$ bits per vertex.\\
\indent For the category of disk graphs, we have $Q=5$, which means this takes $O(n^{0.990839})$ bits per vertex.\\
\\
These are very close to the trivial upper bound of $\left\lceil\frac{n-1}{2}\right\rceil +\log_2(n)$ but are still a small improvement over it when $n$ is large.

\section{Acknowledgements}

Thanks to my supervisor, Oleg Pikhurko, for help and guidance.

Research supported by ERC Grant No. 306493.

\bibliography{report}{}
\bibliographystyle{elsarticle-num}

\end{document}